\documentclass[10pt]{amsart}
\usepackage{amsfonts}
\usepackage{ifthen}
\usepackage{amsthm}
\usepackage{amsmath}
\usepackage{graphicx}
\usepackage{amscd,amssymb,amsthm}
\usepackage{enumerate}
\usepackage{epstopdf}

\newcounter{minutes}
\divide\time by 60
\newcounter{hours}
\multiply\time by 60 \addtocounter{minutes}{-\time}

\title{Tur\'an type inequalities for general Bessel functions}

\author[\'A. Baricz]{\'Arp\'ad Baricz}
\address{\'A. Baricz, Department of Economics,  Babe\c{s}-Bolyai University, Cluj-Napoca 400591, Romania}
\address{Institute of Applied Mathematics, \'Obuda University, 1034 Budapest, Hungary}
\email{bariczocsi@yahoo.com}

\author[S. Ponnusamy]{Saminathan Ponnusamy}
\address{S. Ponnusamy, Indian Statistical Institute, Chennai Centre, Society for Electronic Transactions and Security,
MGR Knowledge City, CIT Campus, Taramani, Chennai 600113, India}
\email{samy@iitm.ac.in, samy@isichennai.res.in}

\author[S. Singh]{Sanjeev Singh}
\address{S. Singh, Department of Mathematics, Indian Institute of Technology Madras, Chennai 600036, India}
\email{sanjeevsinghiitm@gmail.com}

\thanks{$^{\bigstar}$The research of \'A. Baricz was supported by the J\'anos Bolyai Research Scholarship of
the Hungarian Academy of Sciences. The second author is on leave from the Indian Institute of Technology Madras, India. The research of S. Singh was supported by the fellowship of the University Grants Commission, India.}

\newtheorem{theorem}{Theorem}

\begin{document}

\def\thefootnote{}
\footnotetext{ \texttt{File:~\jobname .tex,
          printed: \number\year-0\number\month-0\number\day,
          \thehours.\ifnum\theminutes<10{0}\fi\theminutes}
} \makeatletter\def\thefootnote{\@arabic\c@footnote}\makeatother



\maketitle


\begin{abstract}
In this paper some Tur\'an type inequalities for the general Bessel function, monotonicity and bounds for its logarithmic derivative are derived. Moreover we find the series representation and the relative extrema of the Tur\'anian of general Bessel functions. The key tools in the proofs are the recurrence relations together with some asymptotic relations for Bessel functions.
\end{abstract}

\section{\bf Introduction and the Main Results}
\setcounter{equation}{0}

The Tur\'an type inequalities for orthogonal polynomials and special functions have been studied extensively in the last 70 years. More often, these orthogonal polynomials and special functions are solutions of some second order differential equations. The log-concave/log-convex nature of orthogonal polynomials and special functions have attracted many researchers, and the topic seems to be interesting still nowadays. Some of the results on modified Bessel functions of the first and second kind have been used recently in different problems of applied mathematics and this motivated new researches in this topic. See for example the paper \cite{expo} and the references therein for more details. In this paper we focus on general Bessel functions, called sometimes as cylinder functions. The main motivation to write this paper emerges from the fact about the Bessel functions of the first kind $J_{\nu}$, Bessel functions of the second kind $Y_{\nu}$ and the zeros $c_{\nu,n}$ of general Bessel function satisfying some Tur\'an type inequalities (see \cite{bp1,joshi,lakshmana,lorch,szasz}), and it is natural to ask whether the general Bessel functions $C_{\nu},$ defined by $C_{\nu}(x)=(\cos\alpha) J_{\nu}(x)-(\sin\alpha)Y_{\nu}(x),$ $0\leq\alpha<\pi,$ has some similar properties. As we can see below, from the point of view of Tur\'an type inequalities, the general Bessel function $C_{\nu}$ behaves like $J_{\nu}$ and $Y_{\nu}.$ The results presented in this paper complement the picture about the Tur\'an type inequalities for Bessel functions of the first and second kind. The case $\alpha=0$ corresponds to the case of Bessel function $J_{\nu},$ while $\alpha={\pi}/{2}$ corresponds to the case of Bessel function $Y_{\nu}.$ See \cite{bp1,joshi,lakshmana,szasz} for more details.

\begin{theorem}\label{th1}
The following assertions are valid:
\begin{enumerate}

\item[\bf a.] If $\nu>0$, $0<\alpha<\pi$ and $x\geq c_{\nu,1}$ where $c_{\nu,1}$ is the first positive zero of the general Bessel function $C_{\nu},$ then the following Tur\'an type inequality holds
\begin{equation}\label{turan1}
\Delta_{\nu}(x)=C^2_{\nu}(x)-C_{\nu-1}(x)C_{\nu+1}(x)> \frac{1}{\nu+1}C^2_{\nu}(x).
\end{equation}
Moreover, for $\alpha=0$ the above Tur\'an type inequality holds true for all $x>0$ and $\nu>0$.

\item[\bf b.] If $\nu>1$, $0<\alpha<\pi$ and $x_{\nu}\in (0,c_{\nu,1})$ is the unique root of the equation
    $$C^2_{\nu}(x)-C_{\nu-1}(x)C_{\nu+1}(x)= \frac{1}{\nu+1}C^2_{\nu}(x),$$
then the  Tur\'an type inequality \eqref{turan1} holds true for all $x>x_{\nu}$. Moreover, the inequality \eqref{turan1} is reversed for $0<x<x_{\nu}$.

\item[\bf c.] The function $x\mapsto {xC'_{\nu}(x)}/{C_{\nu}(x)}$ is strictly decreasing on $(c_{\nu,1},\infty)\setminus\Xi$ for all $\nu>0$, $0<\alpha<\pi$, where $\Xi=\left\{c_{\nu,n}\right\}_{n\geq1}$, and $c_{\nu,n}$ denote the $n$th positive zeros of the general Bessel functions $C_{\nu}$. Moreover, if $\alpha=0$ then $x\mapsto {xC'_{\nu}(x)}/{C_{\nu}(x)}$ is strictly decreasing on $(0,\infty)\setminus\Xi$ for all $\nu>0$ and if $0< \alpha<\pi$ then $x\mapsto {xC'_{\nu}(x)}/{C_{\nu}(x)}$ is strictly decreasing on $(x_{\nu},\infty)\setminus\Xi$ for all $\nu>1.$ Furthermore, the following inequality holds true for $\nu>0$, $0<\alpha<\pi$ and $x\in (c_{\nu,1},\infty)\setminus\Xi$
\begin{equation}\label{ine2}
\left[\frac{xC'_{\nu}(x)}{C_{\nu}(x)}\right]^2 >\nu^2-\frac{\nu}{\nu+1}x^2.
\end{equation}
If $\alpha=0,$ then the inequality \eqref{ine2} is valid for all $\nu>0$, $x\in (0,\infty)\setminus\Xi.$ However, if $0< \alpha<\pi,$ then the inequality \eqref{ine2} is valid for all $\nu>1$, $x\in (x_{\nu},\infty)\setminus\Xi,$ and for $x\in (0,x_{\nu})$ it is reversed. The following inequality is also valid for $\nu>1,$ $0<\alpha<\pi$ and $x\in (x_{\nu},c_{\nu,1})$
\begin{equation}\label{bound3}
\frac{xC'_{\nu}(x)}{C_{\nu}(x)} <-\sqrt{\nu^2-\frac{\nu}{\nu+1}x_{\nu}^2}.
\end{equation}
\end{enumerate}
\end{theorem}
The proof of Theorem \ref{th1} will be presented in Section \ref{proof}.

Now we let $\mu=\frac{\nu}{\nu+1}$ and denote by $j_{\nu,n}$, the $n$th positive zero of the Bessel function $J_{\nu}.$ We would like to take the opportunity to mention that by using the particular case of \eqref{ine2} when $\alpha=0,$ it can be shown that for $\nu>0$, $x\in (0,\sqrt{\nu(\nu+1)})$ such that $x\neq j_{\nu-1,n},$ $n\in\mathbb{N}$, we have
$$\frac{J_{\nu}(x)}{J_{\nu-1}(x)}<\frac{\nu-\sqrt{\nu^2-\mu x^2}}{\mu x}<\frac{\nu+\sqrt{\nu^2-\mu x^2}}{\mu x},$$
and this inequality corrects the inequality \cite[eq. 2.20]{bp1}
$$\frac{J_{\nu}(x)}{J_{\nu-1}(x)}\geq \frac{\nu+\sqrt{\nu^2-\mu x^2}}{\mu x},$$
where $\nu>0$, $x\in (0,\sqrt{\nu(\nu+1)})$ such that $x\neq j_{\nu-1,n}$, $n\in\mathbb{N}$.

We also note that the monotonicity of $x\mapsto {xC'_{\nu}(x)}/{C_{\nu}(x)}$ has been proved already by Spigler \cite{spigler} (as it is mentioned in the paper of Elbert and Siafarikas \cite{elbsia}), but only for the intervals $(c_{\nu,n},c_{\nu,n+1})$, $n\in\mathbb{N}.$ Our proof is completely different, which is based on Tur\'an type inequalities and we prove the above monotonicity property for $x\in (x_{\nu},c_{\nu,1})$ and also for $x\in (c_{\nu,n},c_{\nu,n+1})$, $n\in\mathbb{N}$.

We recall the inequality from part {\bf c} of Theorem \ref{th1}
$$
\left[\frac{xC'_{\nu}(x)}{C_{\nu}(x)}\right]^2 <\nu^2-\frac{\nu}{\nu+1}x^2,
$$
where $\nu>1$, $0<\alpha<\pi$ and $x\in(0,x_{\nu})$. We may rewrite the above inequality as
\begin{equation}\label{bound1}
-\sqrt{\nu^2-\frac{\nu}{\nu+1}x^2}<\frac{xC'_{\nu}(x)}{C_{\nu}(x)}<\sqrt{\nu^2-\frac{\nu}{\nu+1}x^2},
\end{equation}
where $\nu>1,$ $0<\alpha<\pi$ and $x\in(0,x_{\nu}).$ We would like to mention that for $\nu>1$, $0<\alpha<\pi$ and $x\in(0,x_{\nu})$ the right-hand side of \eqref{bound1} is better than the earlier known inequality of Laforgia \cite[p. 76]{laforgia}
\begin{equation}\label{bound2}
\frac{xC'_{\nu}(x)}{C_{\nu}(x)}<\nu-\frac{x^2}{2(\nu+1)},
\end{equation}
which is valid for $\nu>0$, $0\leq\alpha<\pi$ and $x\in(0,c_{\nu,1})$. This can be verified by comparing the corresponding right-hand side parts of the last two inequalities. We also note that for $\nu>1$ and $0<\alpha<\pi$ we have $x_{\nu}<\sqrt{\nu(\nu+1)},$ since the expression in the square root in \eqref{bound3} is positive.

Again, it is worth to mention that the inequality \eqref{bound3} is better than the inequality \eqref{bound2} for all $\nu>1$, $0<\alpha<\pi$ and $x\in(x_{\nu},\min\{\sqrt{\nu(\nu+1)},c_{\nu,1}\})$ as the right-hand side of \eqref{bound3} is negative and the right-hand side of \eqref{bound2} is positive on $(x_{\nu},\min\{\sqrt{\nu(\nu+1)},c_{\nu,1}\})$.

The next theorem is about the series representation of the Tur\'anian of general Bessel functions. Clearly, this result provides an alternative proof of the inequality \eqref{turan1}.
\begin{theorem}\label{new}
For $0\leq\alpha<\pi$, $\nu>0$ and $x>c_{\nu,1}$, the following identity holds
\begin{equation}\label{identity}
C^2_{\nu}(x)-C_{\nu-1}(x)C_{\nu+1}(x)= \frac{1}{\nu+1}C^2_{\nu}(x)+2\nu\sum_{i\geq1}\frac{C^2_{\nu+i}(x)}{(\nu+i)^2-1}.
\end{equation}
\end{theorem}

\begin{figure}[!ht]
   \centering
       \includegraphics[width=12cm]{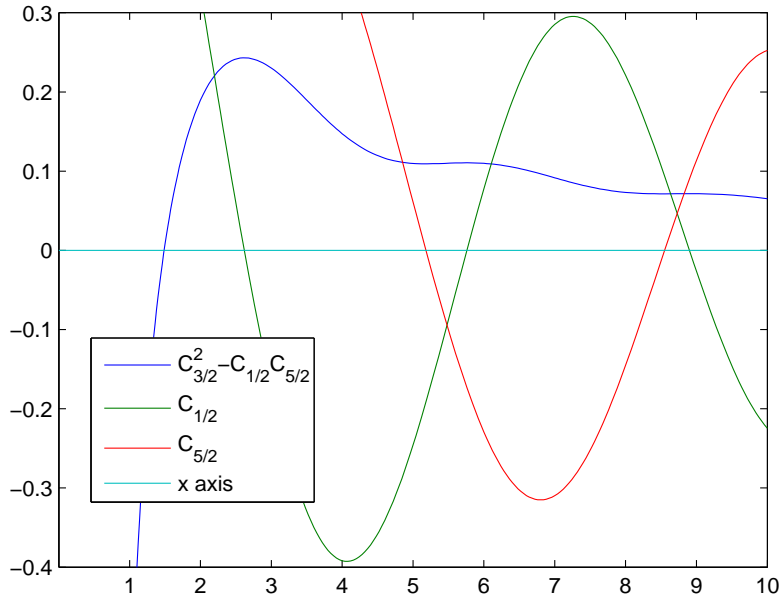}
       \caption{The graph of the functions $\Delta_{3/2},$ $C_{1/2}$ and $C_{5/2}$ for $\alpha=\pi/6$ on $[0,10].$}
       \label{fig1}
\end{figure}

The next result, whose proof will be also presented in Section \ref{proof}, is about the relative extrema of the Tur\'anian of general Bessel functions and is a generalization of the main result from \cite{lakshmana}. Figure \ref{fig1} illustrates this result for $\alpha=\pi/6$ and $\nu=3/2$.

\begin{theorem}\label{th2}
For $0\leq\alpha<\pi$ and $\nu>0$, the relative maxima (denoted by $M_{\nu,k}$) of the function $x\mapsto \Delta_{\nu}(x)$ occurs at the zeros of the function $C_{\nu-1}(x)$ and the  relative minima (denoted by $m_{\nu,k}$) occurs at the zeros of the function $C_{\nu+1}(x)$. Since the values of $M_{\nu,k}$ and $m_{\nu,k}$ can be expressed as
$M_{\nu,k}=\Delta_{\nu}(c_{\nu-1,k})=C^2_{\nu}(c_{\nu-1,k})>0$ and $m_{\nu,k}=\Delta_{\nu}(c_{\nu+1,k})=C^2_{\nu}(c_{\nu+1,k})>0,$ respectively, for $\nu>0$ and $x\geq c_{\nu-1,1}$ the following Tur\'an type inequality is valid:
$$ C^2_{\nu}(x)-C_{\nu-1}(x)C_{\nu+1}(x)>0.$$
\end{theorem}

We also mention that by using the above Tur\'an type inequality in Theorem \ref{th2} for $\alpha=0$ it can be shown that for $x\in (0,\nu)$ such that $x\neq j_{\nu-1,k}$, $k\in\mathbb{N}$, we get
$$\frac{J_{\nu}(x)}{J_{\nu-1}(x)}<\frac{\nu-\sqrt{\nu^2- x^2}}{x}<\frac{\nu+\sqrt{\nu^2- x^2}}{x}.$$
We note that the above inequality corrects the recent known inequality \cite[eq 2.16]{bp1}
$$\frac{J_{\nu}(x)}{J_{\nu-1}(x)}\geq \frac{\nu+\sqrt{\nu^2- x^2}}{x},$$
where $\nu>0$, $x\in (0,\sqrt{\nu(\nu+1)})$ such that $x\neq j_{\nu-1,n}$, $n\in\mathbb{N}.$

Finally, it is worth to mention that the Tur\'anian $\Delta_{\alpha}(x)=C_{\nu,\alpha}^2(x)-C_{\nu,\alpha-1}(x)C_{\nu,\alpha+1}(x),$ where as above $C_{\nu,\alpha}(x)=(\cos\alpha) J_{\nu}(x)-(\sin\alpha)Y_{\nu}(x),$ $0\leq\alpha<\pi,$ is in fact independent of $\alpha.$ Namely, by using some elementary trigonometric identities it can be shown that $\Delta_{\alpha}(x)=(\sin^2 1)\left(J_{\nu}^2(x)+Y_{\nu}^2(x)\right),$ which is clearly strictly positive for all real $\nu$ and $x.$

\section{\bf Proofs of the Main Results}\label{proof}
\setcounter{equation}{0}

\begin{proof}[\bf Proof of Theorem \ref{th1}]
{\bf a.} We first recall the recurrence relation and the derivative formula for general Bessel functions \cite[p. 222]{nist}, which will be used in the sequel
\begin{equation}\label{rec1}
C_{\nu-1}(x)+C_{\nu+1}(x)=\frac{2\nu}{x}C_{\nu}(x),
\end{equation}
and
\begin{equation}\label{diff1}
\frac{d}{dx}(x^{-\nu}C_{\nu}(x))=-x^{-\nu}C_{\nu+1}(x).
\end{equation}
Let us define the normalized general Bessel function by $\Phi_{\nu}(x)=2^{\nu}x^{-\nu}\Gamma(\nu+1)C_{\nu}(x),$ where $\nu>-1$ and $x>0$. Since $C_{\nu}(x)$ is the solution of Bessel differential equation
$$x^2y''(x)+xy'(x)+(x^2-\nu^2)y(x)=0,$$
we see that $\Phi_{\nu}(x)$ satisfies the differential equation
\begin{equation}\label{de2}
x^2\Phi''_{\nu}(x)+(2\nu+1)x\Phi'_{\nu}(x)+x^2\Phi_{\nu}(x)=0.
\end{equation}
Now, if we consider the Tur\'anian $\Theta_{\nu}(x)=\Phi^2_{\nu}(x)-\Phi_{\nu-1}(x)\Phi_{\nu+1}(x),$ then the Tur\'an type inequality \eqref{turan1} is equivalent to $\Theta_{\nu}(x)> 0$. Therefore, it is enough to prove that $\Theta_{\nu}(x)>0.$ Taking into account \eqref{rec1} we have
\begin{equation}\label{rec2}
\frac{x^2\Phi_{\nu+1}(x)}{4\nu(\nu+1)}=\Phi_{\nu}(x)-\Phi_{\nu-1}(x),
\end{equation}
and consequently, in view of \eqref{diff1},
\begin{equation}\label{diff2}
\Phi'_{\nu}(x)=-\frac{x\Phi_{\nu+1}(x)}{2(\nu+1)}=\frac{2\nu}{x}(\Phi_{\nu-1}(x)-\Phi_{\nu}(x)).
\end{equation}
By using the left-hand side of \eqref{diff2}
for $\nu-1$ instead of $\nu$ and the right-hand side of \eqref{diff2}
for $\nu+1$ instead of $\nu$, we obtain that
$$\frac{d}{dx}(x^{2\nu+2}\Theta_{\nu}(x))=-\frac{2}{\nu}x^{2\nu+2}\Phi_{\nu}(x)\Phi'_{\nu}(x).$$
From the above expression, at the roots of $\Phi_{\nu}(x)=0$ we have
$$\frac{d^2}{dx^2}(x^{2\nu+2}\Theta_{\nu}(x))=-\frac{2}{\nu}x^{2\nu+2}\left(\Phi'_{\nu}(x)\right)^2<0$$
and at the roots of $\Phi'_{\nu}(x)=0$, by using \eqref{de2} we obtain
$$\frac{d^2}{dx^2}(x^{2\nu+2}\Theta_{\nu}(x))=\frac{2}{\nu}x^{2\nu+2}\left(\Phi_{\nu}(x)\right)^2>0.$$
These two inequalities show that the relative extrema of $x\mapsto x^{2\nu+2}\Theta_{\nu}(x)$ occurs at the roots of $\Phi_{\nu}(x)=0$ and $\Phi'_{\nu}(x)=0$, respectively. At the roots of $\Phi_{\nu}(x)=0$, we have
$$\Theta_{\nu}(x)=-\Phi_{\nu-1}(x)\Phi_{\nu+1}=\frac{x^2}{4\nu(\nu+1)}\Phi^2_{\nu+1}(x)>0,$$
by using \eqref{rec2} and in the view of \eqref{diff2} at the roots of $\Phi'_{\nu}(x)=0$ one has $\Theta_{\nu}(x)=\Phi^2_{\nu}(x)>0.$ On the other hand, $\Phi_{\nu}(x)=0$ if and only if $C_{\nu}(x)=0$ and in view of \eqref{diff2} we have that
$$\Phi'_{\nu}(x)=0 \Longleftrightarrow \Phi_{\nu+1}(x)=0 \Longleftrightarrow C_{\nu+1}(x)=0.$$
Therefore, the first relative extrema of $x\mapsto x^{2\nu+2}\Theta_{\nu}(x)$ occurs at $x=c_{\nu,1}$, as $c_{\nu,1}<c_{\nu+1,1}$ (since $\nu \mapsto c_{\nu,k}$ is increasing function of $\nu$ \cite[p. 508]{watson}). Since $x\mapsto x^{2\nu+2}\Theta_{\nu}(x)$ has all its relative extrema positive and hence $x^{2\nu+2}\Theta_{\nu}(x)>0$ for all $x\geq c_{\nu,1}$ and $\nu>0,$ which implies that $\Theta_{\nu}(x)>0$ and consequently the Tur\'an type inequality \eqref{turan1} follows for all $\nu>0$ and $x\geq c_{\nu,1}$. Since the first relative extrema of $x\mapsto x^{2\nu+2}\Theta_{\nu}(x)$ occurs at $x=c_{\nu,1},$ which is the point of relative maxima, we conclude that $x\mapsto x^{2\nu+2}\Theta_{\nu}(x)$ is strictly increasing on $(0,c_{\nu,1}).$ Now, if we take $\alpha=0$, then $C_{\nu}(x)=J_{\nu}(x)$ and using the fact that $J_{\nu}(0)=0$ for $\nu>0$ we have $\lim_{x\rightarrow 0^+}x^{2\nu+2}\Theta_{\nu}(x)=0.$ Hence $x^{2\nu+2}\Theta_{\nu}(x)>0$ on $(0,c_{\nu,1})$ and consequently, $\Theta_{\nu}(x)>0$. Therefore in this case, the Tur\'an type inequality \eqref{turan1} holds true for all $\nu>0$ and $x>0.$

{\bf b.} Next consider $0<\alpha<\pi$. Since
$$x^{2\nu+2}\Theta_{\nu}(x)=2^{2\nu}\Gamma(\nu)\Gamma(\nu+2)\left[x^2\left(C^2_{\nu}(x)-C_{\nu-1}(x)C_{\nu+1}(x)\right)-\frac{x^2}{\nu+1}C^2_{\nu}(x)\right],$$
in view of the fact that $J_{\nu}(0)=0$ for $\nu>0$, the asymptotic formula which is valid for $\nu>0$ fixed and $x\rightarrow 0$ \cite[p. 223]{nist}
$$
Y_{\nu}(x)\sim -\frac{1}{\pi}\Gamma(\nu)\left(\frac{x}{2}\right)^{-\nu},
$$
and the limit (see \cite[p. 316]{bp1})
$$\lim_{x\rightarrow 0^+}x^2\left(Y^2_{\nu}(x)-Y_{\nu-1}(x)Y_{\nu+1}(x)\right)=-\infty,$$
where $\nu>1$ is fixed, we obtain that
$$\lim_{x\rightarrow 0^+}x^{2\nu+2}\Theta_{\nu}(x)=-\infty.$$ Therefore, in view of the fact that $x^{2\nu+2}\Theta_{\nu}(x)$ is positive at $x=c_{\nu,1},$ $x^{2\nu+2}\Theta_{\nu}(x)$ tends to $-\infty$ as $x\rightarrow 0^+$ and $x\mapsto x^{2\nu+2}\Theta_{\nu}(x)$ is strictly increasing on $(0,c_{\nu,1})$, we obtain that there exists an unique $x_{\nu}\in (0,c_{\nu,1})$ such that
$$
\begin{cases}
~   x^{2\nu+2}\Theta_{\nu}(x) <0 &   \mbox{for $x\in (0,x_{\nu})$,}\\
~  x^{2\nu+2}\Theta_{\nu}(x) =0 &   \mbox{for $x= x_{\nu}$,} \\
~ x^{2\nu+2}\Theta_{\nu}(x) >0         & \mbox{for $x\in (x_{\nu},c_{\nu,1})$}.
\end{cases}
$$
Hence $\Theta_{\nu}(x)>0$ for $x\in (x_{\nu},c_{\nu,1})$. Consequently, in view of part {\bf a} of this theorem, the Tur\'an type inequality \eqref{turan1} is indeed true for all $\nu>1$ and $x>x_{\nu}$. We also note that $\Theta_{\nu}(x)<0$ for $x\in (0,x_{\nu})$ and hence in this case the inequality \eqref{turan1} is reversed.

{\bf c.} In view of the recurrence relations \cite[p. 222]{nist}
\begin{equation}\label{rec3}C'_{\nu}(x)=C_{\nu-1}(x)-\frac{\nu}{x}C_{\nu}(x)~~ \mbox{and}~~ C'_{\nu}(x)=-C_{\nu+1}(x)+\frac{\nu}{x}C_{\nu}(x),\end{equation}
the Tur\'an expression $\Delta_{\nu}(x)$ can be written as
\begin{equation}\label{turanian}
\Delta_{\nu}(x)=C^2_{\nu}(x)-C_{\nu-1}(x)C_{\nu+1}(x)
=\left(1-\frac{\nu^2}{x^2}\right)C^2_{\nu}(x)+[C'_{\nu}(x)]^2.
\end{equation}
Now, by using this and the Bessel differential equation we get
$$\left(\frac{xC'_{\nu}(x)}{C_{\nu}(x)}\right)'=-\frac{\Delta_{\nu}(x)}{xC^2_{\nu}(x)}.$$
Thus, in view of part {\bf a} and {\bf b} of Theorem \ref{th1}, the monotonicity of $x\mapsto {xC'_{\nu}(x)}/{C_{\nu}(x)}$ follows. We note that the inequality \eqref{ine2} follows from \eqref{turan1} by using the above recurrence relations in \eqref{rec3}.

Since $x\mapsto {xC'_{\nu}(x)}/{C_{\nu}(x)}$ is strictly decreasing on $(x_{\nu},c_{\nu,1})$ for all $\nu>1$, we see that
\begin{equation}\label{ine3}
\frac{xC'_{\nu}(x)}{C_{\nu}(x)}<\lim_{x\rightarrow x_{\nu}}\frac{xC'_{\nu}(x)}{C_{\nu}(x)}.
\end{equation}
Using the relation \eqref{turanian} we obtain that
$$
\left[\frac{xC'_{\nu}(x)}{C_{\nu}(x)}\right]^2=\frac{x^2\Delta_{\nu}(x)}{C^2_{\nu}(x)}+(\nu^2-x^2),
$$
which in view of the fact that (see \cite[p. 78]{laforgia})
$$
C_{\nu}(x)>0 ~~ \mbox{and}~~ C'_{\nu}(x)<0 ~~ \mbox{for}~~ 0<x<c_{\nu,1}
$$
implies that
$$
\frac{xC'_{\nu}(x)}{C_{\nu}(x)}=-\sqrt{\frac{x^2\Delta_{\nu}(x)}{C^2_{\nu}(x)}+(\nu^2-x^2)}.
$$
Taking the limit $x\rightarrow x_{\nu}$ in the above equation and using \eqref{ine3}, we get the inequality \eqref{bound3}.

\end{proof}

\begin{proof}[\bf Proof of Theorem \ref{new}]
Let us recall one of the main results of Landau \cite{landau}: the magnitude of the general Bessel function of order $\nu$ is decreasing in $\nu$ at all of its stationary points. This means that if $\nu>0$ and $x>c_{\nu,1},$ then we have that $|C_{\nu}(x)|<|C_0(x)|<|\tau|,$ where $\tau=C_0(x_1)$ and $x_1$ is the abscissa of the first minimum point of $C_0(x).$ Hence in view of the recurrence relation \eqref{rec1}, all conditions \cite[Theorem 3]{ross} are satisfied and consequently we obtain the identity \eqref{identity}.
\end{proof}

\begin{proof}[\bf Proof of Theorem \ref{th2}]
Using the recurrence relations in \eqref{rec3}, we have
\begin{equation}\label{diff3}
\Delta'_{\nu}(x)=\frac{2}{x}C_{\nu-1}(x)C_{\nu+1}(x).
\end{equation}
Hence the relative extrema of $x\mapsto \Delta_{\nu}(x)$ occurs at the zeros of $C_{\nu-1}(x)$ and $C_{\nu+1}(x)$. From \eqref{diff3}, by using the second recurrence relation in \eqref{rec3} for $\nu-1$ instead of $\nu,$ and \eqref{rec1}, we get
$$\Delta''_{\nu}(x)\big|_{x=c_{\nu-1,k}}=-\frac{4\nu}{c^2_{\nu-1,k}}C^2_{\nu}(c_{\nu-1,k})<0.$$
Similarly, by using the first recurrence relation in \eqref{rec3} for $\nu+1$ instead of $\nu,$ and \eqref{rec1}, we get
$$\Delta''_{\nu}(x)\big|_{x=c_{\nu+1,k}}=\frac{4\nu}{c^2_{\nu+1,k}}C^2_{\nu}(c_{\nu+1,k})>0.$$
The desired conclusion follows.
\end{proof}

\end{document}